%% file: SmartGridComm_longv2.tex
\documentclass[10pt,conference,onecolumn]{IEEEtran}
\usepackage{verbatim}
\usepackage{amssymb}
\usepackage[cmex10]{amsmath}
\usepackage{graphicx}
\usepackage{epsfig}
\usepackage{cite}
\usepackage{psfrag}
\usepackage{comment}
\usepackage{bm}
\usepackage{bbm}
\usepackage{url}
\usepackage{amsthm}
\usepackage{subfig}
\input defs.tex
\input RLDRCdef.tex

\begin{document}
\title{Risk Limiting Dispatch with Ramping Constraints}
\author{\IEEEauthorblockN{Junjie Qin}
\IEEEauthorblockA{ICME\\
Stanford University\\
Stanford, CA 94305\\
Email: jqin@stanford.edu}
\and
\IEEEauthorblockN{Baosen Zhang}
\IEEEauthorblockA{EECS\\
University of California, Berkeley\\
Berkeley, CA 94720\\
Email:zhangbao@eecs.berkeley.edu\\}
\and
\IEEEauthorblockN{Ram Rajagopal}
\IEEEauthorblockA{CEE\\
Stanford University\\
Stanford, CA 94305\\
Email: ramr@stanford.edu}}

\markboth{}%
{Shell \MakeLowercase{\textit{et al.}}: Bare Demo of IEEEtran.cls for Journals}
\maketitle

\begin{abstract}
The increased penetration of renewable resources poses a challenge to reliable operation of power systems. In particular, operators face the risk of not scheduling enough traditional generators in the times when renewable generation becomes lower than expected. This paper studies the optimal trade-off between system risk and the cost of scheduling reserve generators. The problem is modeled as a multi-period stochastic control problem, and in particular explicitly considering the ramping constraints on the generators. The structure of the optimal dispatch control is identified.  Dispatch is efficiently computed utilizing two methods:  (i) solving a surrogate chance constrained program, (ii) a MPC-type look ahead controller.  Utilizing real system data, the chance constrained dispatch is shown to outperform the MPC controller and to be robust to changes in the probability distribution of the  uncertainties involved in the problem.
\end{abstract}


\IEEEpeerreviewmaketitle

\section{Introduction}

Renewable resources are starting to play increasingly prominent roles in electric systems around the world. Major efforts are under way to integrate new renewable resources with existing infrastructure, and this has proven to be a non-trivial task.  The main challenge is in the difference between the time-scale of significant changes in renewable generation output and the time-scale at which traditional generators can ramp up or down.  For example, independent system operators (ISOs) typically schedule generators in hourly blocks, and have limited ramping capabilities. On the other hand, renewable generation can change significantly on intra-hour time-scales, and under high penetration levels, the system may become ramp-limited and incapable of adapting to the real time renewable output.

Uncertainties have an asymmetric effect on system operations in power networks \cite{VWB2011}. For example, if the predicted wind power output is lower than the actual realization, the system may not have enough ramp down capability, and excessive wind needs to be curtailed, resulting in wasted energy and increased cost. Thus the uncertainty has an effect on the operation \emph{cost} of the system.  On the other hand, if the predicted wind power output is  higher than the actual realization, then the system may not have enough ramp up capability and some load needs to be shed. This is a much worse outcome than curtailing generation, and we call this the operation \emph{risk} in the system. Presently, system operators are predominately concerned with the risk, and would purchase enough reserve energy in advance so that the risk that demand cannot be satisfied is negligible.

The current strategy is viable because uncertainties are fairly small (on the order of $1- 2 \%$ of total demand). However, as the penetration of renewables increases, the strategy of procuring enough energy for all possible scenarios is becoming prohibitively expensive \cite{RBDEWV2012}. For example, consider an day ahead of scheduling generators in a system with wind power.  Since there is always some possibility that there would be no wind in the next day, to operate at zero risk, the operator would essentially schedule as if there was no wind. The resulting operation process is neither economical nor sustainable as the required reserves incur in significant emissions. Instead, the operator could accept a tolerable level in system risk obtaining in return significant reduction in required reserves.  One contribution of this paper is to address the problem of how to optimal schedule generation considering the trade-off between operation risk and cost in a ramp-constrained system.

\begin{figure}[ht]
 \centering
 \includegraphics[scale=1]{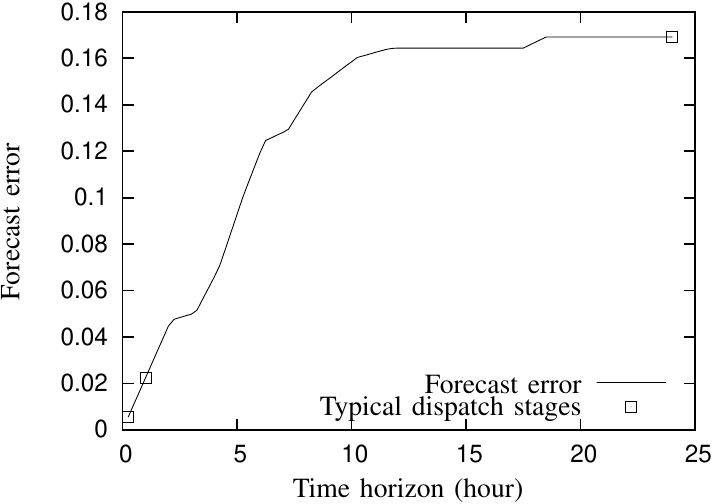}
 \caption{Percentage forecast error v.s. forecast horizon. Data from Iberdrola Renewables.}
 \label{fig:error_update}
 \end{figure}

The second contribution of this paper is to model ramping constrained dispatch as a multi-stage stochastic control problem that attempts to minimize the expected cost of procuring energy, constrained to a desired probability of loss of load and explicit ramping constraints in the generation. The sequential nature of forecast error updates along time is explicitly modeled to accurately capture the information and risk tradeoff. This results in successive dispatch with ramp constraints. Figure \ref{fig:error_update} shows a typical error curve for wind forecasting \cite{RBVW2011}.
As expected, the forecast error decreases as the forecast horizon shortens. We can interpret the improvement as errors being successively revealed as the system moves closer to real time. Therefore as each a decision is made, the operator should take advantage of the new information and perform a recourse action. This type of problem has been studied in inventory control as the newsvendor (or newsboy)  problem \cite{S2012}. However, an important difference is that in power systems a decision is made after the renewable realization is observed for the current period, where as in inventory control a decision is made before the randomness is known. As we will show, the solution is far more complicated when the decision can be based on the current period realization.

Prior work mostly consider independent and identically distributed errors over time, or relaxes the ramping constraints. A two-stage version of the dispatch problem was proposed in \cite{VWB2011} and has been extended to include storage \cite{QSR2012} and network effects \cite{RTZ2012}. Multi-stage stage formulations have mainly been studied under the context of price uncertainty (see e.g. \cite{FRD2012, KP2011}) and sequential reserve calculations \cite{RBV2011}.This paper focuses on the inter-temporal behavior introduced by ramping constraints on the generators. In existing ISO operations, this is the most frequently limiting constraint.  The paper is organized as follows. Sec. \ref{sec:model} describes the problem formulation. Sec. \ref{sec:struc} develops the optimal control structure and practical algorithms.  Sec. \ref{sec:simulation} shows the performance of the dispatch algorithms and Sec \ref{sec:con} concludes the paper.

\section{Model and Problem Formulation} \label{sec:model}
We consider a multi-period stochastic dispatch problem. For simplicity and due to space constraints, we assume that the network can be approximated as a single bus. We note the approach in \cite{RTZ2012} can be used to include the network effects. Let there be $T$ total time periods, and at each time $t=0,1,\dots,T-1$, a dispatch decision $g_t$ is made.  An additional terminal stage, \ie, the $T$th stage, is included for notational convenience. At each stage $t$, the information available to the operator includes (i) all the past dispatches, (ii) all the past realization of load $(l_s:0 \le s \le t)$ and wind power generation $(w_s: 0\le s\le t)$, and (iii) a forecast of the future load $(\hat{l}_{t,s}: t < s \le T)$ and wind power generation $(\hat{w}_{t,s}: t < s \le T)$. Assuming wind power generation is taken at zero cost, we model the wind as negative load, and consider net demand defined as $d_t = l_t - w_t$.

\subsection{Statistical Model}
An important feature of our approach is to model the information update process. In multi-period dispatch problems, the operator obtains better forecasts about the wind and load at fixed future time points as the forecast horizons become shorter. To incorporate this notion, we consider the forecast update model of the form
\begin{equation}
    d_t = \df_{s,t} + \sum_{\tau=s}^{t-1} \tilde{\e}_{\tau,t},
\end{equation}
where the forecast vector $\df_{s}\in \real^{T+1}$  available at time $s$ is defined as
\[
    \df_{s,t} = \begin{cases} d_t & \mbox{if } t \le s,\\
                            \hat{d}_{s,t} & \mbox{if }  t > s,
                \end{cases}
\]
with $\hat{d}_{s,t}$ being the net demand at time $t$ forecasted at stage $s$, and   $\tilde{\e}_\tau\in \real^{T+1}$ is a vector containing {\it marginal forecast errors}, \ie, the forecast error of net demand $d_t$ of the forecast made at time $s$ is
\begin{align*}
    \epsilon_{s,t} = \epsilon_{s+1,t} + \tilde{\e}_{s,t}
                = &\epsilon_{s+2,t}+ \tilde{\e}_{s+1,t}+ \tilde{\e}_{s,t}= \dots\\
                & = \epsilon_{t,t}+ \sum_{\tau=s}^{t-1} \tilde{\e}_{\tau,t} = \sum_{\tau=s}^{t-1} \tilde{\e}_{\tau,t}.
\end{align*}
Note $\tilde{\e}_{s,t}=0$ for $t \le s$, for convenience, we work with the reduced error vector $\e_t \in \real^{T-t}$, such that $\tilde{\e}_{t}^\T = [\zeros^\T, \e_t^\T]^\T$.
It follows that the forecast vector $\df_t$ is updated according to
\[
    \df_{t+1} = \df_t + \C_t \e_t,
\]
where
\[
    \C_t = \begin{bmatrix}
    \zeros_{(t+1)\times (T-t)} \\\eye_{(T-t) \times (T-t)}
    \end{bmatrix}\in \real^{(T+1)\times (T-t)}
\]
is a zero-one matrix that ensures only coordinates of $\df_t$ that corresponding to future periods will be updated. Note a subtle but important difference between our model and a standard inventory model (see, e.g. \cite{BertDPBook}) is that we allow the operator to observe the current error ($e_t$) before the current dispatch ($g_t$) is made.

Empirical studies (\eg, \cite{Makarov2009}) suggests that the forecast error for wind power generation, which is the major source of the uncertainty, follows a (truncated) Gaussian distribution. Here we assume $\e_t$, for each $t$, is zero mean normal random vector with prescribed covariance $\Sigma_t$. In numerical example, we also validate our approach against forecast errors that are not Gaussian.

\subsection{Problem Formulation}
For each unit of conventional generation dispatched, a constant cost $c$ is incurred. We want to control the probability of the error event of excess  demand compared to available generation, i.e. $\{d_t > g_t\}$. The constraint in the system are the ramping constraints of the form
\begin{equation}\label{eqn:ramp}
    \rl \le \g_t - \g_{t-1} \le \ru,
\end{equation}
where $\rl <0$ and $\ru >0$ are constants representing the maximum ramping down and ramping up rate of generation. For convenience, denote the feasible set as
\[
\Gcal(\g)=\{\g'\ge 0: -\rl\le \g'-\g\le\ru\},
\]
where dispatched generation is naturally required to be positive.
The main optimization problem is
\begin{subequations}\label{ccp}
\begin{align}
\mbox{minimize} &\quad \expec\left[\sum_{t=0}^{T-1}\c \g_t + \psi(d_t,g_t)\right]\\
\mbox{subject to} & \quad \df_{t+1} = \df_t + \C_t \e_t, \\
& \quad g_t \in \mathcal G(g_{t-1}) \label{cons_ramp} \\
& \quad g_t = f_t(\df_0,\dots,\df_t,g_0,\dots,g_{t-1}) \label{ccCon5}
\end{align}
\end{subequations}
where $\psi(d_t,g_t)$ is a risk measure that is convex in $g_t$. In particular, we consider LOLP, which takes the form
\[
    \psi (d_t, \g _t ) = \begin{cases} \infty & \mbox{if } \prob (d_t >  g_t)  > \beta_0, \\
                                      0      & \mbox{otherwise,}
                                      \end{cases}
\]
where $\beta_0$ is the tolerance for the chance that the demand is more than the scheduled generation. The penalty function is also mathematically equivalent to the constraint
\begin{equation}
\prob(d_{t}-\g_t \le 0) \ge 1-\beta_0. \label{con1}
\end{equation}
This constraint is convex in current scaler form. In more general cases when net demands at multiple buses or periods are involved, it is convex if the underlying probability distribution is log-concave (e.g. Gaussian, Laplace and see \cite{Boyd04}).
In practice, we  set $\beta_0$ according to the probability that reserves would be needed.
The constraint \eqref{ccCon5} states that the dispatch decision is causal, since $g_t$ is restricted to be a function of all the information that is available up to the current period. Since $g_{t-1}$ and $\df_t$ completely capture the state of the system at time $t$, it turns out that $g_t$ only depend on $g_{t-1}$, $\df_t$, and the statistics of the future errors.\footnote{This constraint is equivalent to requiring the function of control to be adapted to the sigma algebra generated by the current information set. We also assume  the statistics of future errors is available to the controller.}

\section{Dispatch Solutions} \label{sec:struc}
In this section, we study the optimal solution to the dynamic program \eqref{ccp}. The results only depend on the convexity of $\psi$ (so various type of convex risk measures can be used).  Define the cost-to-go function as
\[
J_t(\df_t, \g_{t-1}) = \min_{\g_t\in  \Gcal(\g_{t-1})}Q_t(\df_t,\g_t),
\]
where $Q_t(\df_t,\g_t) $ is the state-action $Q$ function
\[
Q_t(\df_t,\g_t) =  \c \g_t + \psi(\df_{t,t},g_t)+\expec\left[ J_{t+1}(\df_t + \C_t \e_t, \g_{t})\right],
\]
with $J_{T}(\df_T, \g_{T-1})=0$. Regarding the property of functions $J_t$ and $Q_t$, the following convexity result follows from standard arguments in \cite{BertDPBook}.
\begin{proposition}
For time periods $t=0,\dots, T$, the $Q$ function $Q_t(\df_t,\g_t)$ and  cost-to-go function
$J_t(\df_t, \g_{t-1}) $  are both convex in their arguments.
\end{proposition}

Let an unconstrained minimizer of the $Q$ function with respect to action be $S_t(\df_t) \in \argmin_{ \g_{t}\in \real} Q_t(\df_t,\g_t),$
then we have the following structural result about the optimal policy.
\begin{theorem}\label{thm:s}
The threshold rule
\begin{align*}
    \mu^\star_t &(\df_t, \g_{t-1})  \\
    &=\begin{cases}
    S_t(\df_t) & \mbox{if } (\g_{t-1}-\rl)^+ \le S_t(\df_t) \le \g_{t-1} + \ru, \\
    (\g_{t-1}-\rl)^+ & \mbox{if } S_t(\df_t) < (\g_{t-1}-\rl)^+, \\
    \g_{t-1} + \ru & \mbox{if } S_t(\df_t) > \g_{t-1} + \ru,
    \end{cases}
\end{align*}
is an optimal policy to problem \eqref{ccp}.
\end{theorem}

\begin{remark}
The threshold rule in Theorem~\ref{thm:s} follows from the convexity of the sequences of functions $J_t$ and $Q_t$, and is proved by induction similar to threshold rules for celebrated inventory control models. \cite{BertDPBook}.
\end{remark}

The problem of identifying an optimal target function $S_t(\df_t)$ is a hard one since the dimension of the search space, \ie, all functions of $\df_t$, is infinite. As it is common practice in stochastic control literature, by restricting to the class of functions that are linear in $\df_t$, the problem is reduced to a finite dimensional optimization which are usually easier to solve. However, as demonstrated in appendix, the optimization problem for coefficients of the linear function involves minimizing an expectation which cannot be evaluated without sampling.
Other standard approximation approach is based on discretizing the state space (\eg, the potential values $\df_t$ can take) and using a Bellman recursion to compute the cost-to-go value for each state \cite{BertDPBook}. Since $\df_t$ is a multiple dimension vector, discretization suffers from the curse of dimensionality. Furthermore, this approach is entirely numerical, in the sense that other than reporting the target values, the approach gives no intuition and/or insight why such a target value is optimal. It is also difficult to implement in existing dispatch systems as it requires a significant departure from the current practice.

In the rest of this section, we present two methods of approximating $S_t$.
\subsection{Chance Constrained Linear Dispatch} \label{sec:affine}
We replace the problem in \eqref{ccp} with a surrogate chance constrained problem\footnote{For more information on chance constraint programming see, e.g. \cite{KW1994}} as
\begin{subequations}\label{eqn:chance}
\begin{align}
\mbox{minimize} &\quad \expec\left[\sum_{t=0}^{T-1}\c \g_t\right]\\
\mbox{subject to} & \quad \df_{t+1} = \df_t + \C_t \e_t, \\
& \quad \prob(\df_{t,t}-\g_t \le 0) \ge 1-\beta_0, \label{ccCons1}\\
& \quad \prob(\g_t \ge 0)\ge 1-\beta_1, \label{ccCons2}\\
& \quad \prob(\g_t-\g_{t-1}  \ge -\rl)\ge 1-\beta_2, \label{ccCons3}\\
& \quad \prob(\g_t-\g_{t-1}  \le \ru)\ge 1-\beta_3,\label{ccCons4}\\
& \quad g_t = f_t(\df_0,\dots,\df_t,g_0,\dots,g_{t-1})
\end{align}
\end{subequations}
Comparing \eqref{ccp} and \eqref{eqn:chance}, we replaced inequalities in the hard constraint \eqref{cons_ramp} by probabilistic constraints \eqref{ccCons2}, \eqref{ccCons3}, and \eqref{ccCons4}, respectively. There are two motivations for using a chance constrained surrogate: i) the constraint \eqref{con1} is already in the form of a chance constraint, and without the ramping constraints Eqs.~\eqref{eqn:chance} and \eqref{ccp} are equivalent  if we let $\beta_k \to 0$ for $k=1,2,3$ \cite{RBVW2011}; ii) we can derive a linear dispatch solution to \eqref{eqn:chance} via a second order cone program.

Let
$\df = \begin{bmatrix} \df_0^\T & \dots & \df_T^\T \end{bmatrix}^\T, $
and
$ \e = \begin{bmatrix} \e_0^\T & \dots & \e_{T-1}^\T\end{bmatrix}^\T.$
Then the forecast update can be summarized as
\[
\df = \A \df_0 + \C \e,
\]
where
$    \A = \begin{bmatrix}
    \eye_{(T+1)\times (T+1)}&
    \dots&
    \eye_{(T+1)\times (T+1)}
    \end{bmatrix}^\T, $
\[
    \C = \begin{bmatrix}
    \zeros_{(T+1)\times n_0} & \zeros_{(T+1)\times n_1} & \dots & \zeros_{(T+1)\times n_{T-1}}\\
    \C_0 & \zeros_{(T+1)\times n_1} & \dots & \zeros_{(T+1)\times n_{T-1}}\\
    \C_0 & \C_1 & \dots & \zeros_{(T+1)\times n_{T-1}}\\
    \vdots & \vdots & \ddots & \vdots\\
    \C_0 & \C_1 & \dots & \C_{T-1}
    \end{bmatrix},
\]
with $n_t = T-t$.
Affine policies are of the form
\begin{equation}
g_t = \sum_{\tau=0}^{t-1} \G_{t,\tau} \e_{\tau} + a_t,
\end{equation}
where $\G_{t,\tau}$ and $a_t$ are parameters to be decided.
Let $\gb = \begin{bmatrix} \gb_0^\T & \dots & \gb_{T-1}^\T \end{bmatrix}^\T$, then    $\gb = \G \e + \a$,
where
$ \a = \begin{bmatrix} a_0 & a_1 &\dots& a_{T-1}\end{bmatrix}^\T$, and
\[
\G = \begin{bmatrix}
\zeros_{1\times n_0} & \zeros_{1\times n_1} & \dots & \zeros_{1\times n_{T-2}} &\zeros_{1\times n_{T-1}} \\
\G_{1,0} & \zeros_{1\times n_1} & \dots & \zeros_{1\times n_{T-2}} &\zeros_{1\times n_{T-1}} \\
\G_{2,0} & \G_{2,1} & \dots &  \zeros_{1\times n_{T-2}} &\zeros_{1\times n_{T-1}} \\
\vdots & \vdots & \ddots & \vdots & \vdots\\
\G_{T-1,0} & \G_{T-1,1} & \dots & \G_{T-1,T-2} & \zeros_{1\times n_{T-1}}
\end{bmatrix}, 
\]
The objective function can be written as
$
\expec\left[\sum_{t=0}^{T-1}\c \g_t\right] = \c \ones^\T \a.
$
To summarize the linear inequality constraints (before taking the probability), let
\[
    \H^d_1 = \begin{bmatrix}
    \deB_0^\T &&&& \zeros_{1\times (T+1)}\\
    & \deB_1^\T &&&\zeros_{1\times (T+1)}\\
    && \ddots &&\vdots\\
    &&& \deB_{T-1}^\T&\zeros_{1\times (T+1)}
    \end{bmatrix},
\]
where $\delta_i\in \real^{1\times(T+1)}$ is the $i$th coordinate vector, and
\[
    \H^g_3 = \begin{bmatrix}
    1 & -1 \\
    & \ddots & \ddots \\
    && 1 & -1
    \end{bmatrix}, \quad 
    \H^g_4 = - \H^g_3.
\]

Then the feasible region for the original program can be written as
\[
    \H^d \df + \H^g \gb - \y \le \zeros,
\]
where
\[
    \H^d = \begin{bmatrix}
    \H^d_1 \\
     \zeros_{T\times (T+1)^2}\\
     \zeros_{(T-1)\times (T+1)^2}\\
 \zeros_{(T-1)\times (T+1^2)}
    \end{bmatrix}, 
    \H^g = \begin{bmatrix}
    -\eye_{T\times T}\\
    -\eye_{T\times T}\\
    \H^g_3\\
    \H^g_4
    \end{bmatrix},
    \y = \begin{bmatrix}
    \zeros_{T\times 1}\\
    \zeros_{T\times 1}\\
    \rl \ones_{(T-1)\times 1}\\
    \ru \ones_{(T-1)\times 1}
    \end{bmatrix}.
\]

Plugging in expressions for $\df$ and $\gb$ yields
\[
    \H^d \df + \H^g \gb - \y  = \h + \P \e \le \zeros,
\]
where
\[
    \h = \H^d \A \df_0 - \y +\H^g \a,\quad
    \P = \H^d \C + \H^g \G.
\]
Let $\P_i^\T$ be the $i$th row of $\P$, then $\P_i^\T \e$ is a Gaussian random variable with zero mean and standard deviation $\|\Sb^{\frac{1}{2}} \P_i\|_2$. Here $\Sb$ is the covariance matrix of $\e$, which can be obtained from the covariance of $\e_t$. The chance constraints then can be expressed as
\[
     h_i + \alpha_i\|\Sb^{\frac{1}{2}} \P_i\|_2 \le 0
\]
where
\[
\alpha_i = \sqrt{2}\,\mathrm{erf}^{-1}(1-2\beta_i'),
\]
with $\beta_i'$ equals the corresponding $\beta_k$, $k=0,1,2,3$, according to the chance one allows constraint $i$ to be violated.
Thus the chance constrained program can be posed as a second order cone program of parameters $\G$ and $\a$
\begin{subequations}\label{ccp_cone}
\begin{align}
\mbox{minimize} &\quad \c \ones^\T \a\\
\mbox{subject to} & \quad h_i + \alpha_i\|\Sb^{\frac{1}{2}} \P_i\|_2 \le 0,
\end{align}
\end{subequations}
and solved efficiently with standard convex optimization solver.

\subsection{Lookahead Policies} \label{sec:lookahead}

In this section we analyze and develop lookahead policies similar to MPC controllers proposed in the literature \cite{VHRW08}, except our error assumptions are distinct  and closed form results are obtained for these controllers. We start by identifying the close-form expressions for $S_t(\df_t)$ for the last two periods of the problem, which in turn will lead to a one-step lookahead policy. In addition to LOLP, we also develop results based on value of lost load (VOLL) penalty function of the form
\[
    \psi(d_t, \g_t) = q(d_t-\g_t)^+,
\]
where $q$ is a positive constant which is typically much larger than $c$. The following results assume that $q>3c$. The proofs of the results in this section are in the appendix.
\begin{lemma}
With LOLP penalty, the optimal target at period $T-2$ is
\begin{align}
&S_{T-2}\label{eqn:SLOLP}\\
&=\max\left\{\df_{T-2,T-2}, \df_{T-2,T-1}-\ru +\Phi_{T-2,1}^{-1}(1-\beta_0) , S'_{T-2}\right\},\nonumber
\end{align}
where $S'_{T-2}$ is the solution to the equation
\begin{align}
& 1+ \ind(\g_{T-2}> \rl)\Phi_{T-2,1}(\g_{T-2}-\rl)\label{onestepLOLP}\\
 &  \quad \quad- \ind(\g_{T-2}>\rl) \df_{T-2,T-1} \phi_{T-2,1}(\g_{T-2}-\rl)=0,\nonumber
\end{align}
where $\Phi_{T-2,1}(\cdot)$ and $\phi_{T-2,1}(\cdot)$ are the cdf and pdf of $\e_{T-2,1}$, respectively.
\end{lemma}
\begin{lemma}
With VOLL penalty, the optimal target at period $T-2$ is
\begin{align}
&S_{T-2}\label{eqn:SVOLL}\\
&=\max\left\{\df_{T-2,T-2}, \df_{T-2,T-1}-\ru +\Phi_{T-2,1}^{-1}\left(\frac{q-2c}{q}\right) , S'_{T-2}\right\},\nonumber
\end{align}
where $S'_{T-2}$ is the solution to the equation
\begin{align}
&(2c-q)  \nonumber\\
&\quad+ c\ind(\g_{T-2}> \rl) \Phi_{T-2,T-1}(\g_{T-2}-\rl-\df_{T-2,T-1})  \label{onestepVOLL}\\
&\quad \qquad+(q-c) \Phi_{T-2,T-1}(\g_{T-2}+\ru-\df_{T-2,T-1})=0.\nonumber
\end{align}
Furthermore,
\begin{equation}\label{app1step}\textstyle{
\tilde{S}_{T-2} = \max\left\{\df_{T-2,T-2}, \df_{T-2,T-1}-\ru+\Phi_{T-2,1}^{-1}\left(\frac{q-2c}{q-c}\right)\right\},
}\end{equation}
is a conservative approximate to $S_{T-2}$ such that
\begin{equation}\label{eq:vollbound}
0\le \tilde{S}_{T-2} - S_{T-2}\le \Phi_{T-2,1}^{-1}\left(\frac{q-2c}{q-c} \right)- \Phi_{T-2,1}^{-1}\left(\frac{q-2c}{q}\right).\end{equation}

\end{lemma}

Note here ``conservative'' means potentially more energy is dispatched reducing the risk of shortfall.
\begin{remark}[One-step lookahead policy]
Computing $S_t$ using Eqn.~\eqref{onestepLOLP} or Eqn.~\eqref{onestepVOLL} gives the optimal one-step lookahead policy for corresponding penalty function. For simplicity, in case of VOLL penalty is used, one may instead use Eqn.~\eqref{app1step} which gives an accurate conservative approximate to the optimal target when $q$ is sufficiently larger than $c$.
\end{remark}

\begin{remark}[Heuristic multi-step generalization of \eqref{app1step}]
Formulae \eqref{app1step} gives a conservative approximate solution to \eqref{onestepVOLL} and has a clear interpretation. Looking one step ahead, the target generation level $S_t$ should be at least $\df_{t,t}$ to meet the demand at current stage, and $\df_{t,t+1}-\ru$ plus some uncertainty margin to be able to meet the demand at next stage. This intuition has a simple multi-step generalization:
    \[
        \tilde{S}_t = \max\left\{\df_{t,t}, \max_{\tau> t} \left\{\tilde{\Phi}_{t,\tau-t}^{-1}\left(\frac{q-2c}{q-c}\right)-(\tau-t)\ru+\df_{t,\tau}\right\}\right\},
    \]
where $\tilde{\Phi}_{t,\tau-t}(\cdot)$ is the cdf of $\epsilon_{t,\tau} = \sum_{m= t}^{\tau-1} \e_{m,\tau-m}$.
\end{remark}

When LOLP penalty is used to compute the dispatch rule, it is clear that at the unfavor event that the dispatched generation cannot supply the net demand, certain cost is incurred for the system operator to maintain the power balance. This cost may be associated with dispatching fast generation like spinning reserve to cover the shortfall. The cost for each unit of such shortfall is precisely the per unit VOLL penalty $q$. This suggests that the LOLP tolerance $\beta_0$ may be selected according to generation cost $c$ and VOLL penalty $q$.
\begin{remark}[Relating LOLP and VOLL based targets] Comparing Eqn.~\eqref{eqn:SLOLP} with Eqn.~\eqref{eqn:SVOLL} and Eqn.~\eqref{app1step} suggests that one step lookahead policy derived from LOLP penalty with $1-\beta_0 = (q-2c)/(q-c)$ is at least as conservative as one step lookahead policy derived from VOLL penalty function.
\end{remark}
\section{Case Studies} \label{sec:simulation}

BPA 2011 \cite{BPAData} data set for wind and load sequences is used for the numerical experiment. The $5$ minute data is aggregated in each hour. We solve the problem for each day under consideration, thus  the decision horizon $T$ for the example is $24$. In the total $365$ days of the year, $100$ days are picked at random, whose net demand profile is depicted in Figure 2(a).

The forecast errors have variances that are increasing with the forecast horizon. The empirical relation between forecast error and forecast horizon is calculated using the empirical curve shown in Figure~\ref{fig:error_update}. Figure 2(b) demonstrates the forecast vector $\df_t$ available to the operator at different times of the day (8 a.m. and 16 p.m., at a particular day of the year). All the previous realizations of the net demand have been observed, thus $\df_t$ records exact net demand realized; the future net demand is forecasted such that the forecast error has increasing variance with forecast horizon, resulting in the 95\% confidence interval increasing over time for each forecast.

\begin{figure}[htbp]
\centering
\subfloat[Net demand from BPA data set]{
\label{1a}
  \includegraphics[width=.4\textwidth]{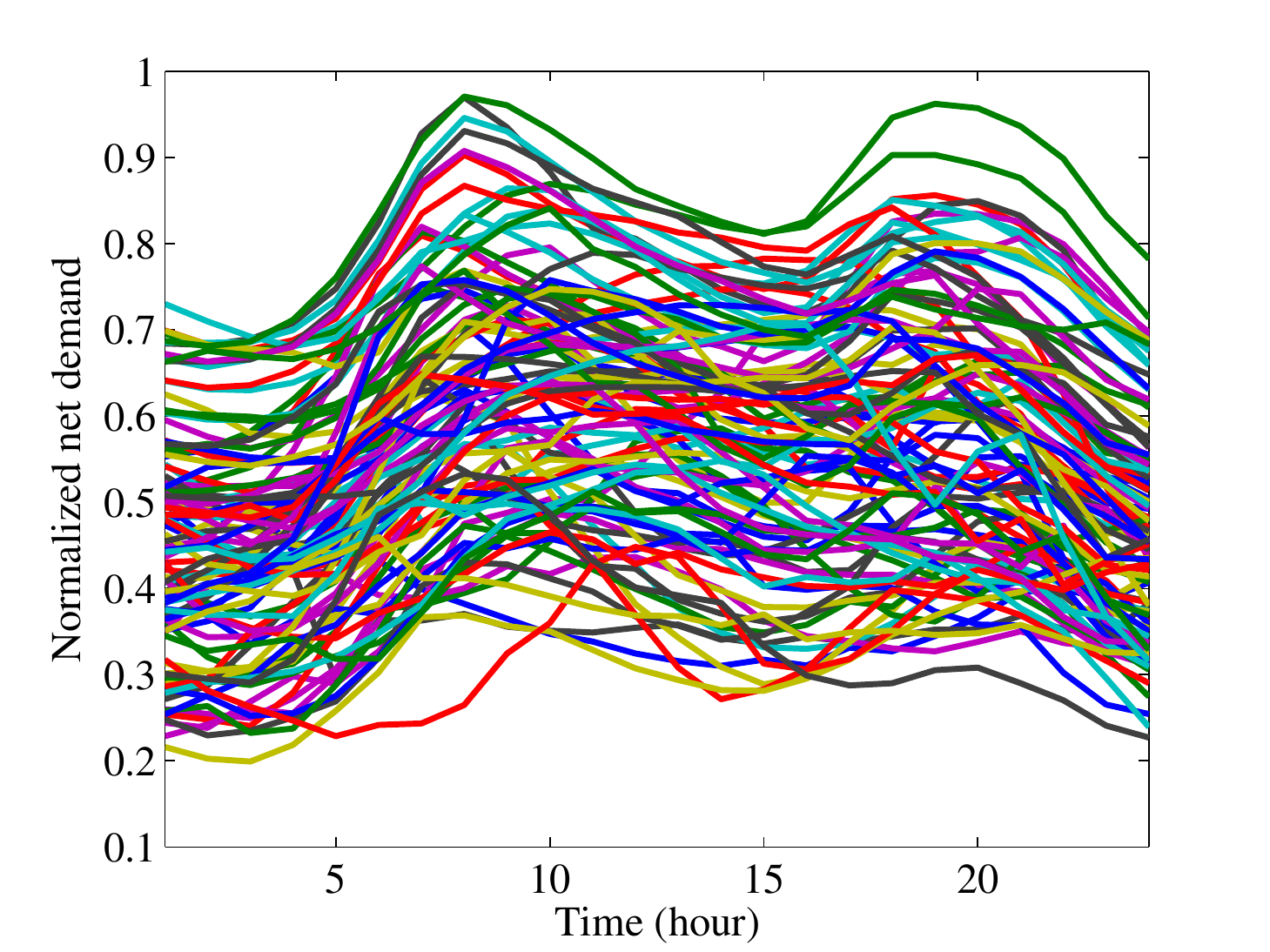}
}\\
\subfloat[Net demand forecast and confidence interval]{
\label{1b}
\includegraphics[width=.4\textwidth]{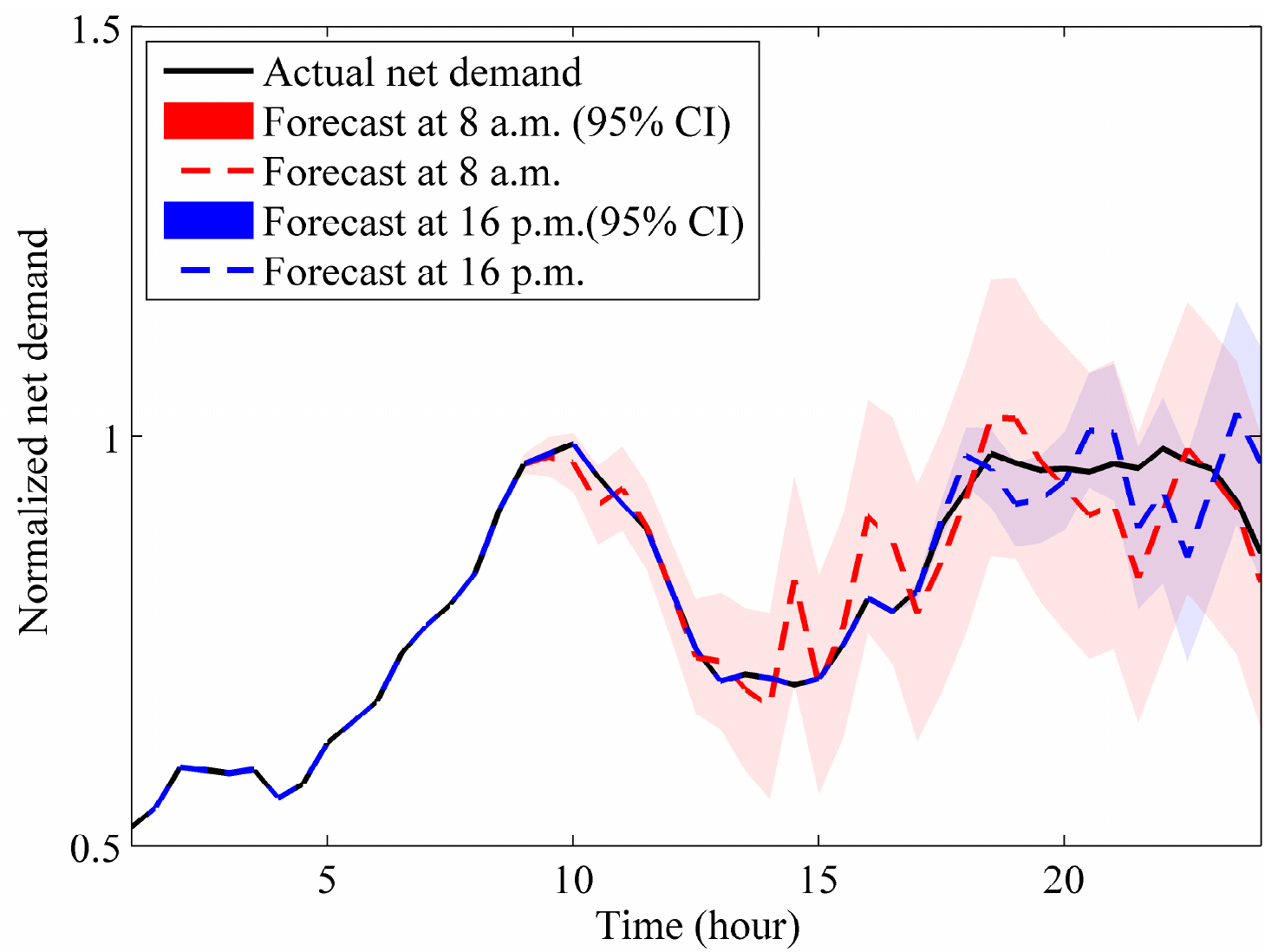}
}
\caption{Net demand data, forecast and forecast error.}
\end{figure}

To make a fair comparison for different policies, we use VOLL penalty to evaluate the performance of control. That is when a loss of load event occurs, the additional cost according to per unit VOLL penalty $q$ is accounted in the cost. Per unit costs  $c=50$, $q=2000$ are set following typical practice in CAISO. Wind generation is scaled such that total wind generation over the day is $p\%$ of the total load, where $p$ is the penetration level desired to simulate a scenario with penetration level $p$.  The ramping rates are set so that $\ru=\rl$ and equal to $4/5$ of the average of the sequence $\{|d_{t+1}-d_t|\}_{t=0}^{T-1}$ to avoid trivial problem instances where the ramping rate is too large (such that the constraint is never binding), or the ramping rate is too small (such that the ramping constraint is always binding).

The performance of different control policies are compared against the oracle cost, which is obtained via solving a deterministic convex optimization using the actual $\mathbf{d}$ sequence. Note that the oracle cost is a lower bound for the cost that can be achieved by any policies that is derived from the partial information contained in the forecast and the distribution of forecast error. In particular, this lower bound does not change while the forecast becomes less accurate, which is likely to be the case when the wind penetration level increases.

\begin{figure}[htbp]
\centering
\includegraphics[width=.4\textwidth, height=5cm]{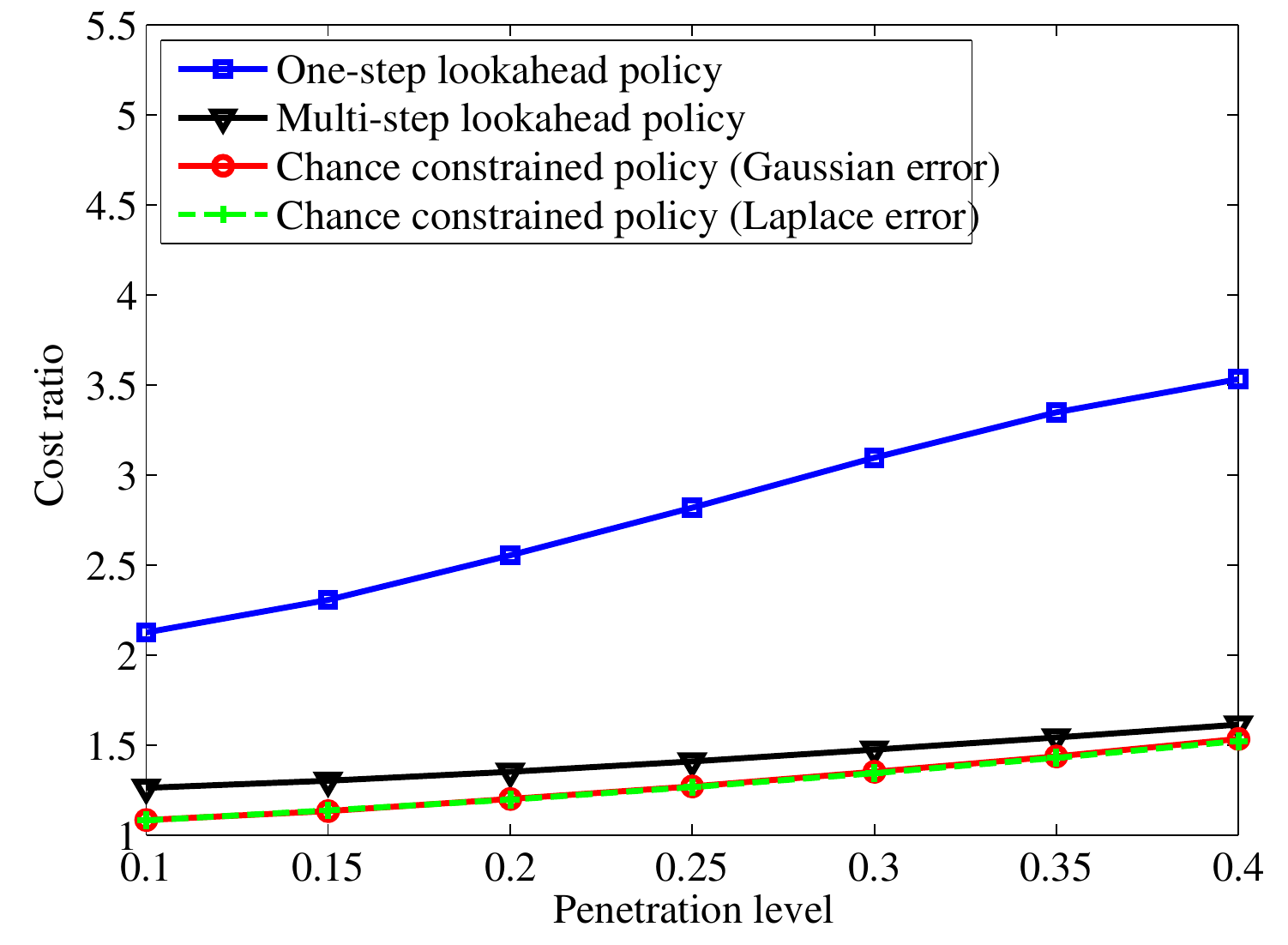}
\caption{Performance of different policies over the test data.}
\label{fig:cost}
\end{figure}

Setting $\beta_k=0.03$ for each $k$, we evaluate the performance of lookahead and chance constrained policies using the metric of cost ratio, \ie, the cost of using the specific policy under consideration divided by the oracle cost. For an efficient policy with mild amount of uncertainty, we expect the cost ratio to be close to $1$. The ratio should  increase slowly when penetration level increases. Figure~\ref{fig:cost} depicts the resulting average cost ratio for various policies. While the chance constrained policy is derived assuming Gaussian forecast error, we also evaluate the cost of the policy using Laplace error (which may be the distribution of the forecast error when certain simple predictor like persistence is used \cite{SEG2011}) of the same standard deviation. The result suggests that chance constrained program works well for the lower penetration range $p\le 0.2$. When the penetration level increases, the cost ratio between chance constrained policy and the oracle cost increases slowly. The cost curves of chance constrained policy against Gaussian and Laplace error almost overlap each other, which indicates the chance constrained policy can be robust against different error distributions. One step lookahead policy (using target \eqref{app1step}\footnote{We have also evaluated the performance of other one step lookahead target formulas. Target derived from LOLP performs worse than VOLL as the cost evaluation is based on VOLL. VOLL targets \eqref{eqn:SVOLL} and \eqref{app1step} produces indifferent cost curves in all our simulations because of the merit of \eqref{eq:vollbound}.}) performs poorly even when the penetration level is low. The multi-step heuristic generalization derived in this paper has a much better performance and somewhat closer to the chance constrained algorithm.

\section{Conclusion and Future Work} \label{sec:con}
This paper formulates and analyzes the problem of risk limiting dispatch with ramping constraints. We model explicitly the forecast update and generation ramping constraints, which are two central elements of dispatch that couple the decision problems for multiple time periods. With dynamic programming, structural results are obtained. Efficient algorithms are then devised to solve the dispatch problem numerically. A case study with real wind data is conducted to illustrate the procedure and the effectiveness of the proposed methods. In future work, we will generalize the framework to incorporate network constraints and additional forward contract markets.

\appendices
\section{Compute linear target with hard constraints}
Using a simple example, we illustrate the complexity of computing linear target with hard constraints in dynamic programming. Suppose we have identified a linear function representing the optimal target for stage $t+1$, denoted as $S_{t+1}(\df_{t+1}) = \elv^\T_{t+1} \df_{t+1}+b_{t+1}$, then the optimal control at stage $t+1$ is
\[
g^\star_{t+1} (\df_{t+1},\g_t)= [\elv^\T_{t+1} \df_{t+1}+b_{t+1}]_{(\g_{t}-\rl)^+}^{\g_t+\ru},
\]
where $[x]_l^u= \min(\max(x,u),l)$. It follows that the state-action $Q$ function at stage $t$ is
\[
Q_t(\df_t, g_{t}) = cg_{t} + \psi(\df_{t,t},g_t) +\expec Q_{t+1}(\df_t+\C_t\e_t,g^\star_{t+1} (\df_{t+1},g_t)),
\]
and then the resulting optimization program for computing $J_t$ involves minimizes $Q_t(\df_t, \elv^\T_{t} \df_{t}+b_t)$ over $\elv_t\in \real^{T+1}$ and $b_t\in\real$. However, it is clear that the second term is a expectation of a function  of a piecewise linear function of the optimization variables, which in general cases cannot be evaluated without sampling.

\section{Lookahead Policies}
\subsection{Derivation of lookahead policies for VOLL penalty}
The cost-to-go function at the last stage is $J_{T}$ is $ J_T(\df_T, \g_{T-1}) =0$. At stage $T-1$, the state-action $Q$ function is
\[
    Q_{T-1}(\df_{T-1},\g_{T-1}) =  q (\df_{T-1,T-1}-\g_{T-1})^++\c \g_{T-1},
\]
and the cost-to-go function is
\[
    J_{T-1}(\df_{T-1}, \g_{T-2}) = \begin{cases}
    \c \df_{T-1,T-1} & \mbox{if } (\g_{T-2}-\rl)^+ <  \df_{T-1,T-1} \le \g_{T-2} + \ru, \\
    \c (\g_{T-2}-\rl)^+ & \mbox{if }  \df_{T-1,T-1} \le (\g_{T-2}-\rl)^+, \\
    q \df_{T-1,T-1} + (\c-q)(\g_{T-2} + \ru) & \mbox{if }  \df_{T-1,T-1} > \g_{T-2} + \ru. \\
    \end{cases}
\]
It follows that the state-action $Q$ function at stage $T-2$ is
\begin{align*}
    &Q_{T-2}(\df_{T-2},\g_{T-2}) \nonumber\\
    =&  q (\df_{T-2,T-2}-\g_{T-2})^++\c \g_{T-2} \nonumber \\
     &+ \expec\left[ J_{T-1}(\df_{T-2} + \C_{T-2} \e_{T-2}, \g_{T-2})\right]\\
     =&  q (\df_{T-2,T-2}-\g_{T-2})^++\c \g_{T-2}  \nonumber\\
     & + \c \df_{T-2,T-1} \prob\left[ (\g_{T-2}-\rl)^+-\df_{T-2,T-1} < \e_{T-2,1}\le \g_{T-2} + \ru-\df_{T-2,T-1}\right] \nonumber\\
    &  +\c (\g_{T-2}-\rl)^+ \prob\left[   \e_{T-2,1}\le (\g_{T-2}-\rl)^+-\df_{T-2,T-1}\right]\nonumber\\
     &  +\left[q \df_{T-2,T-1} + (\c-q)(\g_{T-2} + \ru)\right]\prob\left[  \e_{T-2,1}> \g_{T-2} + \ru-\df_{T-2,T-1}\right]\quad\quad\quad\nonumber \\
     & + \expec \Bigg\{\e_{T-2,1} \bigg\{\c\ind\left[ (\g_{T-2}-\rl)^+-\df_{T-2,T-1} < \e_{T-2,1}\le \g_{T-2} + \ru-\df_{T-2,T-1}\right] \nonumber \\
     &  + q\ind\left[  \e_{T-2,1}> \g_{T-2} + \ru-\df_{T-2,T-1}\right]\bigg\}\Bigg\}.
\end{align*}
A subgradient with respect to $\g_{T-2}$ can be calculated
\begin{align*}
\nabla_{\g} &Q_{T-2}(\df_{T-2},\g_{T-2}) =\\& (2c-q) -q\ind(\df_{T-2,T-2}>\g_{T-2}) + c\ind(\g_{T-2}> \rl) \Phi_{T-2,T-1}(\g_{T-2}-\rl-\df_{T-2,T-1}) +(q-c) \Phi_{T-2,T-1}(\g_{t-2}+\ru-\df_{T-2,T-1}).
\end{align*}
By setting the above expression to $0$, we recover the optimal target $S_{T-2}$ as the root of the corresponding equation.

We proceed to analyze the equation $\nabla_{\g} Q_{T-2}(\df_{T-2},\g_{T-2}) =0$. Put
\[
f(\g_{T-2}) = (2c-q) -q\ind(\df_{T-2,T-2}>\g_{T-2}) + c\ind(\g_{T-2}> \rl) \Phi_{T-2,T-1}(\g_{T-2}-\rl-\df_{T-2,T-1}) +(q-c) \Phi_{T-2,T-1}(\g_{t-2}+\ru-\df_{T-2,T-1}),
\]
it is clear that $f$ is nondecreasing in $\g_{T-2}$. We show $S_{T-2}\ge \df_{T-2,T-2}$, and \[\df_{T-2,T-1}-\ru + \Phi^{-1}_{T-2,1}\left(\frac{q-2c}{q}\right)\le S_{T-2} \le \df_{T-2,T-1}-\ru + \Phi^{-1}_{T-2,1}\left(\frac{q-2c}{q-c}\right).\]
First suppose $S_{T-2}< \df_{T-2,T-2}$, then\footnote{Note the following derivation holds actually for all subgradients.}
\[
f(S_{T-2})  \le  2c-2q+c+(q-c) =2c-q <0,
\]
if $2c<q$. Now suppose $S_{T-2}< \df_{T-2,T-2}$ is in force, we have
\[
f\left(\df_{T-2,T-1}-\ru + \Phi^{-1}_{T-2,1}\left(\frac{q-2c}{q}\right)\right) < (2c-q)+c \frac{q-2c}{q} +(q-c)\frac{q-2c}{q} =0,
\]
and
\[
f\left(\df_{T-2,T-1}-\ru + \Phi^{-1}_{T-2,1}\left(\frac{q-2c}{q-c}\right)\right) >(2c-q) + (q-c) \frac{q-2c}{q-c}=0,
\]
which completes the proof.
\subsection{Derivation of lookahead policies for LOLP penalty}
For notational ease, we work with extended reals $\overline{\real} = \real \cup \{-\infty, +\infty\}$ and write LOLP penalty as
\[
    \psi (d_t, \g _t ) = \begin{cases} M & \mbox{if } \prob (d_t >  g_t)  > \beta, \\
                                      0      & \mbox{otherwise,}
                                      \end{cases}
\]
where $M= \infty$.
At the last stage, $t=T$, the cost-to-go function is $J_T(\df_T, \g_{T-1}) =0$.
At stage $T-1$, the state-action $Q$ function is
\[
    Q_{T-1} (\df_{T-1}, \g_{T-1})  = cg_{T-1} + M \ind(\prob(d_{T-1}> g_{T-1} > \beta)),
\]
where given $d_{T-1}$ is already observed at stage $T-1$, \ie, $d_t=\df_{T-1,T-1}$ is a deterministic quantity conditioning on the current information, the last indicator is equivalent to $\ind (\df_{T-1,T-1}> g_{T-1})$. The ``complication'' involved here will have a clear value when we roll back one stage. The cost-to-go function is then
\begin{align*}
J_{T-1} (\df_{T-1}, g_{T-2})& = \min_{g_{T-1} \in \Gcal(g_{T-2})} c g_{T-1} + M \ind ( \prob(\df_{T-1,T-1> g_{T-1}}) > \beta)\\
& = c \max(\df_{T-1,T-1}, (g_{T-2}-\rl)^+) + M \ind ( \prob(\df_{T-1,T-1} > g_{T-2}+\ru) > \beta).
\end{align*}
The next to last stage, $T-2$, has associated state-action $Q$ function
\begin{align*}
Q_{T-2} (\df_{T-2}, g_{T-2})& = c g_{T-2} + M \ind(\prob(\df_{T-2,T-2} > g_{T-2})>\beta) + \expec  J_{T-1}(\df_{T-2}+C_{T-2} \e_{T-2}, \g_{T-2}) \\
& = c g_{T-2} + M \ind(\prob(\df_{T-2,T-2} > g_{T-2})>\beta) + M \ind ( \prob(\df_{T-2,T-1}+\e_{T-2,1} > g_{T-2}+\ru) > \beta)\\
& \qquad \qquad+ \expec  c \max(\df_{T-2,T-1}+\e_{T-2,1}, (g_{T-2}-\rl)^+) ,
\end{align*}
where we used the fact that
\[
    \expec \ind ( \prob(\df_{T-2,T-1}+\e_{T-2,1} > g_{T-2}+\ru) > \beta) = \ind ( \prob(\df_{T-2,T-1}+\e_{T-2,1} > g_{T-2}+\ru) > \beta)
\]
because the quantity inside the indicator is deterministic. Then the cost-to-go function is
\begin{align*}
J_{T-2} (\df_{T-2}, g_{T-3}) &= \min_{\g_{T-2}\in \Gcal(\g_{T-3})} \Bigg\{ c g_{T-2} + M \ind(\prob(\df_{T-2,T-2} > g_{T-2})>\beta) + M \ind ( \prob(\df_{T-2,T-1}+\e_{T-2,1} > g_{T-2}+\ru) > \beta)\\
& \qquad \qquad\qquad\qquad+ \expec  c \max(\df_{T-2,T-1}+\e_{T-2,1}, (g_{T-2}-\rl)^+)\Bigg\} \\
&= \min_{\overset{\overset{\g_{T-2}\in \Gcal(\g_{T-3})}
                            {g_{T-2}>\df_{T-2,T-2} }}
                            {\prob(\df_{T-2,T-1}+\e_{T-2,1} > g_{T-2}+\ru) > \beta}}
                                                                            \Bigg\{ c g_{T-2} +  \expec  c \max(\df_{T-2,T-1}+\e_{T-2,1} (g_{T-2}-\rl)^+)\Bigg\} \\
&= \min_{\overset{\overset{\g_{T-2}\in \Gcal(\g_{T-3})}
                             {g_{T-2}>\df_{T-2,T-2} }}
                            { g_{T-2}\ge \df_{T-2,T-1}-\ru +\Phi^{-1}_{T-2,1}(1-\beta)    }}
                                                                            \Bigg\{ c g_{T-2} +  \expec  c \max(\df_{T-2,T-1}+\e_{T-2,1} (g_{T-2}-\rl)^+)\Bigg\}.
\end{align*}
Put $f(g) = c g +  \expec  c \max(\df_{T-2,T-1}+\e_{T-2,1}, (g-\rl)^+)$, then
\begin{align*}
    f(g) & = c g + c \int_{\real}  \max(\df_{T-2,T-1}+e (g-\rl)^+) \phi_{T-2,1}(e) de \\
    & = c g + c (g-\rl)^+ \int_{-\infty}^{(g-\rl)^+}   \phi_{T-2,1}(e) de  + c \int_{(g-\rl)^+}^\infty  (\df_{T-2,T-1}+e) \phi_{T-2,1}(e) de
\end{align*}
It follows that
\[
    \grad_g f(g) = c + c \ind(g> \rl)\Phi_{T-2,T-1}(g-\rl)  - c\ind(g>\rl) \df_{T-2,T-1} \phi_{T-2,1}(g-\rl).
\]
Thus $S'_{T-2}$ is the root of $\grad_g f(g)=0$. The additional constraints involving LOLP implies the optimal target at stage $T-2$ is
\[
S_{T-2} = \max\left\{\df_{T-2,T-2}, \df_{T-2,T-1}-\ru +\Phi^{-1}_{T-2,1}(1-\beta) , S'_{T-2}\right\}.
\]

\bibliographystyle{IEEETran}
\bibliography{RLD}

\end{document}

%% file: defs.tex
   \usepackage{epstopdf}
\newcommand{\BEAS}{\begin{eqnarray*}}
\newcommand{\EEAS}{\end{eqnarray*}}
\newcommand{\BEQ}{\begin{equation}}
\newcommand{\EEQ}{\end{equation}}
\newcommand{\BIT}{\begin{itemize}}
\newcommand{\EIT}{\end{itemize}}

\newcommand{\eg}{{\it e.g.}}
\newcommand{\ie}{{\it i.e.}}

\newcommand{\ones}{\mathbf 1}

\newcommand{\ind}{\mathbbm{1}}

\newcommand{\argmin}{\mathop{\rm argmin}}

\newcommand{\real}{\mathbb{R}}

\newcommand{\prob}{\mathbb{P}}
\newcommand{\expec}{\mathbb{E}}

\newcommand{\grad}[1]{\nabla#1} 

\newtheorem{theorem}{Theorem}[section]
\newtheorem{lemma}[theorem]{Lemma}
\newtheorem{proposition}[theorem]{Proposition}

\newtheorem{remark}[theorem]{Remark}

\usepackage{nomencl}

\makenomenclature

\usepackage{makeidx}
\makeindex

\usepackage{tikz}
\usetikzlibrary{shapes,backgrounds,calc}

\makeatletter
\tikzset{circle split part fill/.style  args={#1,#2}{%
 alias=tmp@name, 
  postaction={%
    insert path={
     \pgfextra{%
     \pgfpointdiff{\pgfpointanchor{\pgf@node@name}{center}}%
                  {\pgfpointanchor{\pgf@node@name}{east}}%
     \pgfmathsetmacro\insiderad{\pgf@x}
      \fill[#1] (\pgf@node@name.base) ([xshift=-\pgflinewidth]\pgf@node@name.east) arc
                          (0:180:\insiderad-\pgflinewidth)--cycle;
      \fill[#2] (\pgf@node@name.base) ([xshift=\pgflinewidth]\pgf@node@name.west)  arc
                           (180:360:\insiderad-\pgflinewidth)--cycle;            
         }}}}}
 \makeatother

\def\e{\epsilon}

\usepackage{algorithmic}
\usepackage{graphicx, subfig}
\usepackage{arydshln}

%% file: RLDRCdef.tex
\newcommand{\g}{g}
\newcommand{\T}{\intercal}
\newcommand{\df}{\hat{\mathbf{d}}}

\newcommand{\zeros}{\mathbf{0}}
\newcommand{\C}{\mathbf{C}}
\renewcommand{\c}{c}
\newcommand{\deB}{\bm{\delta}}
\newcommand{\ru}{\overline{r}}
\newcommand{\rl}{\underline{r}}

\newcommand{\Gcal}{\mathcal{G}}
\renewcommand{\e}{\mathbf{e}}
\newcommand{\eye}{\mathbf{I}}

\newcommand{\A}{\mathbf{A}}
\renewcommand{\H}{\mathbf{H}}
\newcommand{\y}{\mathbf{y}}
\newcommand{\G}{\mathbf{G}}
\newcommand{\gb}{\mathbf{g}}
\renewcommand{\a}{\mathbf{a}}
\newcommand{\h}{\mathbf{h}}
\renewcommand{\P}{\mathbf{P}}
\newcommand{\Sb}{\bm{\Sigma}}
\def\elv{\bm{\ell}}

%% file: SmartGridComm_longv2.bbl
\begin{thebibliography}{10}
\providecommand{\url}[1]{#1}
\csname url@samestyle\endcsname
\providecommand{\newblock}{\relax}
\providecommand{\bibinfo}[2]{#2}
\providecommand{\BIBentrySTDinterwordspacing}{\spaceskip=0pt\relax}
\providecommand{\BIBentryALTinterwordstretchfactor}{4}
\providecommand{\BIBentryALTinterwordspacing}{\spaceskip=\fontdimen2\font plus
\BIBentryALTinterwordstretchfactor\fontdimen3\font minus
  \fontdimen4\font\relax}
\providecommand{\BIBforeignlanguage}[2]{{%
\expandafter\ifx\csname l@#1\endcsname\relax
\typeout{** WARNING: IEEEtran.bst: No hyphenation pattern has been}%
\typeout{** loaded for the language `#1'. Using the pattern for}%
\typeout{** the default language instead.}%
\else
\language=\csname l@#1\endcsname
\fi
#2}}
\providecommand{\BIBdecl}{\relax}
\BIBdecl

\bibitem{VWB2011}
P.~Varaiya, F.~Wu, and J.~Bialek, ``Smart operation of smart grid:
  Risk-limiting dispatch,'' \emph{Proceedings of the IEEE}, vol.~99, no.~1, pp.
  40--57, 2011.

\bibitem{RBDEWV2012}
R.~Rajagopal, J.~Bialek, C.~Dent, R.~Entriken, F.~F. Wu, and P.~Varaiya, ``Risk
  limiting dispatch: Empirical study,'' in \emph{12th International Conference
  on Probabilistic Methods Applied to Power Systems}, 2012.

\bibitem{RBVW2011}
R.~Rajagopal, E.~Bitar, F.~F. Wu, and P.~Varaiya, ``{Risk-Limiting Dispatch for
  Integrating Renewable Power},'' \emph{International Journal of Electrical
  Power and Energy Systems, to appear}, 2012.

\bibitem{S2012}
W.~J. Stevenson, \emph{Operations Management}, 11th~ed.\hskip 1em plus 0.5em
  minus 0.4em\relax McGraw and Hill, 2012.

\bibitem{QSR2012}
\BIBentryALTinterwordspacing
J.~Qin, H.~Su, and R.~Rajagopal, ``Risk limiting dispatch with fast ramping
  storage,'' \emph{Submitted to IEEE Transactions on Automatic control}, 2012.
  [Online]. Available: \url{http://arxiv.org/abs/1212.0272}
\BIBentrySTDinterwordspacing

\bibitem{RTZ2012}
R.~Rajagopal, D.~Tse, and B.~Zhang, ``Network risk limiting dispatch: Optimal
  control and price of uncertainty,'' \emph{arXiv:1212.4898}, 2012.

\bibitem{FRD2012}
A.~{Faghih}, M.~{Roozbehani}, and M.~A. {Dahleh}, ``{On the Economic Value and
  Price-Responsiveness of Ramp-Constrained Storage},'' \emph{ArXiv e-prints},
  2012.

\bibitem{KP2011}
J.~H. Kim and W.~B. Powell, ``Optimal energy commitments with storage and
  intermittent supply,'' \emph{Operations Research}, 2011.

\bibitem{RBV2011}
R.~Rajagopal, E.~Bitar, F.~F. Wu, and P.~Varaiya, ``{Risk Limiting Dispatch of
  Wind Power},'' in \emph{{{Proceedings of the American Control Conference
  (ACC)}}}, 2012.

\bibitem{BertDPBook}
D.~P. Bertsekas, \emph{Dynamic Programming and Optimal Control}.\hskip 1em plus
  0.5em minus 0.4em\relax Athena Scientific, 2007.

\bibitem{Makarov2009}
Y.~Makarov, C.~Loutan, J.~Ma, and P.~de~Mello, ``Operational impacts of wind
  generation on california power systems,'' \emph{Power Systems, IEEE
  Transactions on}, vol.~24, no.~2, pp. 1039--1050, 2009.

\bibitem{Boyd04}
S.~Boyd and L.~Vandenberghe, \emph{Convex Optimization}.\hskip 1em plus 0.5em
  minus 0.4em\relax Cambridge Press, 2004.

\bibitem{KW1994}
P.~Kall and S.~W. Wallace, \emph{Stochastic Programming}.\hskip 1em plus 0.5em
  minus 0.4em\relax Wiley, 1994.

\bibitem{VHRW08}
A.~N. Venkat, I.~A. Hiskens, J.~B. Rawlings, and S.~J. Wright, ``Distributed
  mpc strategies with application to power system automatic generation
  control,'' \emph{IEEE Transaction on Control System Technology}, 2008.

\bibitem{BPAData}
\BIBentryALTinterwordspacing
{Bonneville Power Administration}. Wind generation \& total load in the {BPA}
  balancing authority. [Online]. Available:
  \url{http://transmission.bpa.gov/Business/Operations/Wind/}
\BIBentrySTDinterwordspacing

\bibitem{SEG2011}
H.-I. {Su} and A.~{El Gamal}, ``{Limits on the Benefits of Energy Storage for
  Renewable Integration},'' \emph{ArXiv e-prints}, Sep. 2011.

\end{thebibliography}
